\documentclass[10pt,brazil,english]{article}
\usepackage{amsfonts}
\usepackage{infocomp}
\usepackage{times}
\usepackage{amsmath}
\usepackage{amssymb}
\usepackage[T1]{fontenc}
\usepackage{comment}
\usepackage[english, portuguese]{babel}
\addto\captionsportuguese{
\renewcommand{\figurename}{Figure}
\renewcommand{\tablename}{Table}
\renewcommand{\refname}{References}}
\usepackage[utf8]{inputenc}
\usepackage{multirow}
\usepackage{lscape}
\usepackage{rotating}
\usepackage{setspace} 
\usepackage[table,xcdraw]{xcolor}
\usepackage{scalefnt}
\usepackage{graphicx}
\usepackage{hyperref}
\usepackage{subfigure}
\usepackage{enumerate}
\usepackage{caption}
\usepackage[sort,compress]{cite}
\usepackage{cleveref,bm}
\usepackage{multicol,multirow,booktabs}

\usepackage{cleveref,bm}
\usepackage{amssymb}
\usepackage{algorithm}
\usepackage[]{algpseudocode}

\usepackage[alf,
abnt-full-initials=no,
abnt-emphasize=bf,								
abnt-and-type=e,
abnt-url-package=url,
abnt-repeated-author-omit=yes,abnt-etal-list=0]{abntex2cite}

\usepackage{fancyhdr}
\usepackage{mathtools}

\newtheorem{thm}{THEOREM}[section]

\newcommand{\overbar}[1]{\mkern 1.5mu\overline{\mkern-1.5mu#1\mkern-1.5mu}\mkern 1.5mu}

\newcommand{\om}{\omega}

\usepackage{rangecite}

\title{ Two-stage stochastic programming approach for path planning problems under travel time and availability uncertainties \\ 
}

\address{
$^{1}$Department of Industrial and Systems Engineering, Wayne State University\\
$^{2}$Department of Industrial and Systems Engineering, Virginia Tech\\
$^{3}$Ground Vehicle Robotics, U.S. Army CCDC Ground Vehicle Systems Center\\ 
$^{4}$Aptiv, Troy, Michigan\\\\
} 

\author{Saravanan Venkatachalam$^{1}$, Manish Bansal$^{2}$, \\ Jonathon M. Smereka $^{3}$, Joseph Lee$^{4}$ }

\abstract{Significant advances in sensing, robotics, and wireless networks have enabled the collaborative utilization of autonomous aerial, ground and underwater vehicles for various applications. However, to successfully harness the benefits of these unmanned ground vehicles (UGVs) in homeland security operations, it is critical to efficiently solve UGV path planning problem which lies at the heart of these operations. Furthermore, in the real-world applications of UGVs, these operations encounter uncertainties such as incomplete information about the target sites, travel times, and the availability of vehicles, sensors, and fuel. This research paper focuses on developing algebraic-based-modeling framework to enable the successful deployment of a team of vehicles while addressing uncertainties in the distance traveled and the availability of UGVs for the mission.}


\begin{document}
\pagestyle{fancy} 

\maketitle
\newpage

\section{Introduction} 

Intelligence, surveillance and reconnaissance (ISR) are critical missions within military operations, and modern-day combat zones pose important challenges for ISR [\cite{chgral}, \cite{zaloga2011unmanned}, \cite{krishnamoorthy2012uav}]. ISR operations are maintained through effective and efficient information collection. Unmanned ground vehicles (UGVs) are an important asset for ISR, target engagements, convoy operations for resupply missions, search and rescue, environmental mapping, disaster area surveying and mapping. Depending upon the nature of the missions, UGVs are preferred over other collection assets. Some of the instances where UGVs are preferred include: unsuitable terrain for human or unmanned aerial vehicles (UAVs), harsh and hostile environment, tedious information collection process for humans, and many more.

Despite the numerous advantages of UGVs, their size and limited payload capacity lead to fuel constraints and therefore, they are required to make one or more refueling stops in a long mission. Moreover, these operations encounter unknown terrain or obstacles, resulting in uncertainty in the fuel (or time) required to travel among different points of interests (POIs); for example, in a hostile terrain with improvised explosive devices (IEDs), conducting anti-IED sweeps and explosive ordinance disposal can lead to unexpected delays for UGVs. In fact, in many applications, even the locations of the POIs are not precisely known (uncertain) due to inaccurate a-priori map or imperfect and noisy exteroceptive sensory information or perturbations; for example, in a fire monitoring application, the POIs of UGVs change based on the random propagation of the fire \cite{casbeer2005forest}. Likewise, other types of system uncertainties include availability of UGVs with specific attributes such as sensors or terrain-compatible vehicle dynamics.

Due to these challenges, to successfully harness the benefits of the UGVs, it is critical to efficiently solve the UGV path-planning problem (UGVPP). Note that the NP-hard problems such as multiple traveling salesman problem (TSP) and distance-constrained vehicle routing problem are special cases of the UGVPP. In this paper, we consider extensions of UGVVPP with aforementioned uncertainties, and refer to this class of problems as Stochastic AVPP (S-AVPP). Motion planning literature \cite{dadkhah2012survey} for autonomous vehicles (AVs) classifies uncertainties into four categories: vehicle dynamics, knowledge of environment, operational environment, and pose information. Uncertainties in operational environment like wind and atmospheric turbulences suite to UAVs and pose information is regarding localization of UGVs. This paper focuses on UGVPP with uncertainties in vehicle dynamics and knowledge of environment. Some of the previous works include: analysis of robustness of modular vehicle fleet considering uncertainties in demand of vehicles \cite{li2017robustness}; path planning for multiple UGVs for deterministic data using heuristic \cite{bellingham2003multi}; single vehicle path planning problems for UGVs considering environmental uncertainties [\cite{evers2014online}, \cite{evers2014robust}. The path planning problem for multiple UGVs considering uncertainties in vehicle dynamics and environmental uncertainties simultaneously using algebraic modelling framework is new to the literature. Furthermore, these algorithms will also be applicable to tackle similar challenges arising in path-planning for UAVs and underwater vehicles, which are used for crop monitoring, ocean bathymetry, forest fire monitoring, border surveillance, and disaster management.

\begin{figure}[!hbtp]
	\begin{center}
		\includegraphics[scale=0.4]{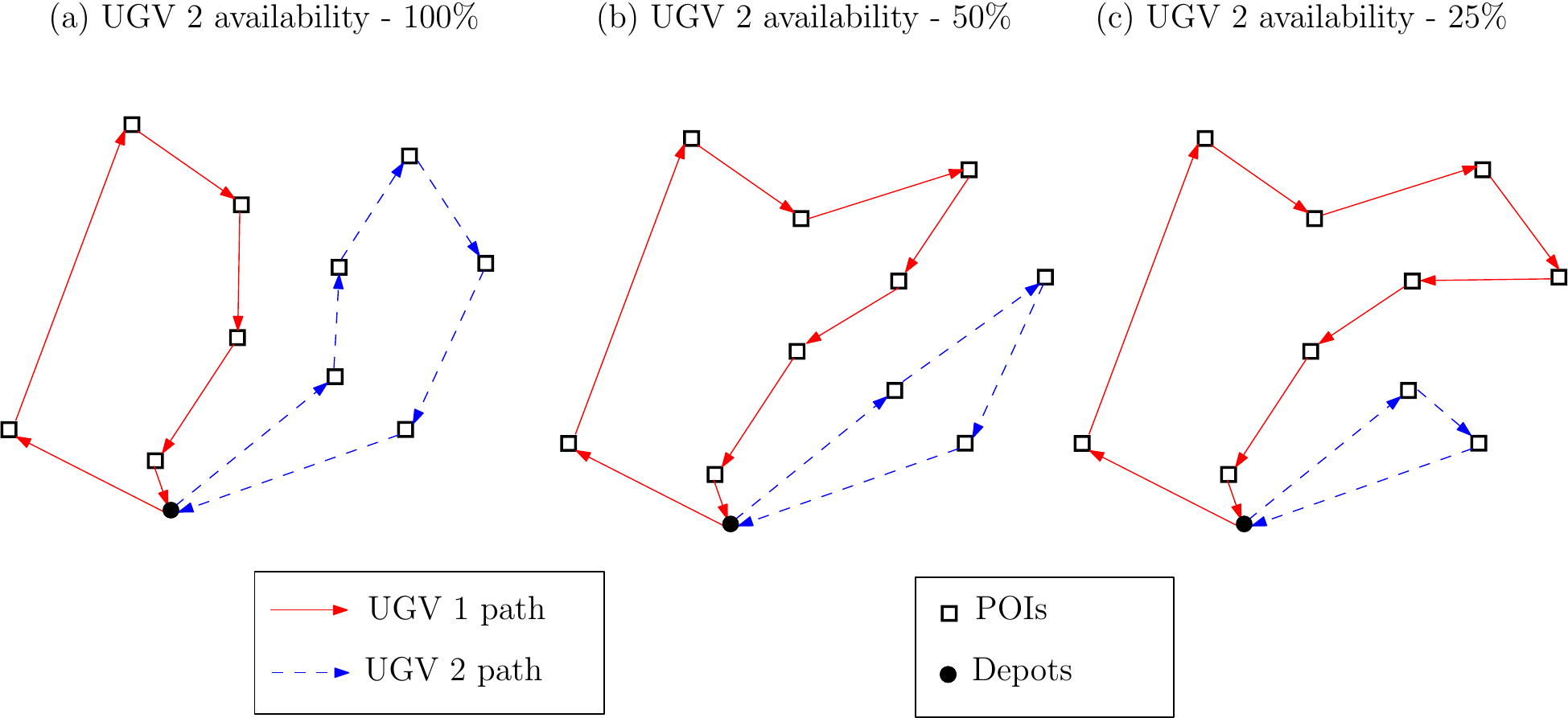}
	\end{center}
	\caption{An illustration of considering the availability of UGVs as uncertain in the UGVPP. (a) Optimal solution for a deterministic UGVPP which is sub-optimal for the UGVs when considering uncertainty for the availability of UGVs. (b-c) Optimal solutions for stochastic UGVPP instances having different chances of availability for UGV2. Note that as the chances of the availability of UGV2 reduces, the number of assigned POIs to UGV2 also reduces.}
	\label{Fig1}
\end{figure}

	\section{Notation} \label{sec:notations}
Let $T= \{t_1,\dots,t_n\}$ denote the set of points of interests (POIs), let $d_0$ denote the depot where a set of heterogeneous unmanned ground vehicles (UGVs) $M:=\{1,\ldots,|M|\}$, each with fuel capacity $F_m, m \in M$, are initially stationed, let $\overbar{D}=\{d_1,\dots, d_k\}$ denote the set of additional $k$ depots or refueling sites, and let $D = \overbar D \bigcup \{d_0\}$. All of the $|M|$ UGVs stationed at the depot $d_0$ are assumed to be fueled to capacity. The model formulations are defined on a directed graph $G = (V, E)$, where $V = T\cup D$ denotes the set of vertices and $E$ denotes the set of edges joining any pair of vertices. We assume that $G$ does not contain any self-loops. For each edge $(i,j) \in E$, we let $c_{ij}$ and $\hat{f}_{ij}^m$ represent the travel cost and the nominal fuel that will be consumed by UGV $m \in M$ while traversing the edge $(i,j)$. We remark that $\hat{f}_{ij}^m\,$ is directly computed using the length of the edge $(i,j)$ and the fuel economy of the UGV. Additional notations that will be used in the mathematical formulation are as follows: for any set $S \subset V$, $\delta^+(S)=\{(i,j) \in E: i\in S, j\notin S\}$ and $\delta^-(S)=\{(i,j) \in E: i\notin S, j\in S\}$. When $S = \{i\}$, we shall simply write $\delta^+(i)$ and $\delta^-(i)$ instead of $\delta^+(\{i\})$ and $\delta^-(\{i\})$, respectively.

The notation introduced next is for describing the uncertainty associated with the UGVs' fuel consumption. Let $\bm f$ denote a discrete random variable vector representing the fuel consumed by any UGV to traverse any edge in $E$. The vector $\bm f$ has $|E|\times |M|$ components, one for each edge, and the random variable in the vector $\bm f$ corresponding to edge $(i,j)$  is denoted by $f_{ij}$. Let $\Omega$ denote the set of scenarios for $\bm f$, where $\omega \in \Omega$ represents a random event or realization of the random variable $\bm f$ with a probability of occurrence $p(\omega)$. We use $f_{ij}^m(\omega)$ to denote the fuel consumed by an UGV $m$ when traversing the edge $(i,j)$, and $\bm f(\omega) = \bigg(\big\{f^1_{ij}(\omega)\big\}_{(i,j)\in E}, \ldots, \big\{f^{|M|}_{ij}(\omega)\big\}_{(i,j)\in E} \bigg)$ to denote the random vector for the realization $\omega \in \Omega$. Finally, we use $\mathbb E$ to denote the expectation operator, i.e., $\mathbb{E}_\Omega(\alpha) = \sum_{\omega \in \Omega}p(\omega)\alpha$. Table \ref{tab:notations} lists all the notations introduced in this section for ease of reading. In the next section, we present two-stage stochastic program formulations using the notation introduced in this section. 

\begin{table}[htbp]
	\centering
	\scalebox{0.6}{
	\begin{tabular}{ll}
		\toprule 
		Symbol & explanation \\ 
		\midrule 
		$T = \{t_1, \dots, t_n\}$ & set of $n$ POIs  \\
		$d_0$ & depot where $m$ UGVs are initially stationed \\ 
		$\overbar{D} = \{d_1, \dots, d_k\}$ & additional depots/refueling sites \\ 
		$D = \overbar D \bigcup \{d_0\}$ & set of depots \\ 
		$F_m$ & fuel capacity of an `$m$' UGV \\ 
		$G = (V, E)$ & directed graph with $V = T\cup D$ \\ 
		$c_{ij}$ & travel cost for the edge $(i,j) \in E$ \\
		$\hat{f}_{ij}^m$ & nominal fuel consumed by the UGV to traverse the edge $(i,j) \in E$ \\
		$\bm f$ & random variable vector representing the fuel consumed by any UGV\\
		$p_j^m$ & profit or incentive for visiting a POI by an UGV $m$\\        
		$\Omega$ & set of scenarios for $\bm f$ \\
		$\omega \in \Omega$ & realization of random variable $\bm f$ \\
		$p(\omega)$ & probability of occurrence of $\omega$ \\ 
		$\mathbb E$ & expectation operator \\
		\bottomrule
	\end{tabular}}
	\caption{Table of notations}
	\label{tab:notations}
\end{table}

\section{Mathematical formulation} \label{sec:formulation}
The first-stage decision variables represent `here-and-now' decisions that are determined before the realization of randomness, and second-stage decisions are determined after scenarios representing the uncertainties are presented. The first-stage decision variables in the stochastic program are used to compute the initial set of routes for each of the UGVs such that either each POI is visited by only one of the UGVs or all the UGVs collect maximum incentives from the POIs, while ensuring that no UGV ever runs out of fuel as it traverses its route. The fuel constraint for each UGV in the first-stage is enforced using the nominal fuel consumption value $\hat f_{ij}^m$ for each edge $(i,j) \in E$. For a realization $\omega \in \Omega$, the second-stage decision variables are used to compute the recourse costs  that must be added to the first-stage routes based on the realized values of $f_{ij}^m(\omega)$ for all $(i,j) \in E$ and $m \in M$. 

Specifically, the first-stage decision variables are as follows: each edge $(i,j)\in E$ is associated with a variable $x_{ij}^m$ that equals $1$ if the edge $(i,j)$ is traversed by `$m$' UGV, and $0$ otherwise. We let $\bm x \in \{0,1\}^{|E| \times |M|}$ denote the vector of all decision variables $x_{ij}^m$. There is also a flow variable $z_{ij}$ associated with each edge $(i,j) \in E$ that denotes the total nominal fuel consumed by any UGV as it starts from depot $i$ and reaches the vertex $j$. Additionally, for any $A \subseteq E$, we let $x^m(A) = \sum_{(i,j)\in A} x_{ij}^m$. Analogous to the variable $x_{ij}^m$ in the first stage, we define a binary variable $y_{ij}^m(\omega)$ for each edge $(i,j)\in E$. The variables $y_{ij}^m(\omega)$ are used to define the refueling trips needed for any vehicle when the route defined by the first-stage feasible solution $\bm x$ is not feasible for the realization $\omega \in \Omega$. 

\section{Formulation 1} \label{sec:for2} 

\textit{Given a team of heterogeneous UGVs (each UGV with a different capacity and travel time between POIs) and multiple refueling depots, a set of target POIs to visit and stochastic travel times or fuel consumption, find a path for each UGV such that each POI site is visited by at most one UGV, and the overall distance traveled by the UGVs is minimized}. 

\subsection{Objective function} \label{subsec:obj2} The objective function for the two-stage stochastic programming model is the sum of the first-stage travel cost and the expected second-stage recourse cost. The second-stage recourse cost for a realization $\omega \in \Omega$ of the fuel consumption of the vehicles is the cost of the additional refueling trips that are required for the realization $\omega$. The recourse cost is a function of the first-stage routing decision $\bm x$ and the realization $\omega$. Letting the recourse cost be denoted by $\beta(\bm x, \bm f(\omega))$, the objective function for the two-stage stochastic optimization problem is given by:

\begin{flalign}
\min \,\, C, \text{ where } C \triangleq \sum\limits_{\substack{(i,j) \in E \\ m \in M}} \hat f_{ij}^m x_{ij}^m + \mathbb{E}_{\Omega} \left[ \beta(\bm x,\bm f) \right]  \nonumber \\ =  \sum\limits_{\substack{(i,j) \in E \\ m \in M}} \hat f_{ij}^m x_{ij}^m + \sum_{\omega \in \Omega} p(\omega) \beta(\bm x, \bm f(\omega)). \label{eq:obj2}
\end{flalign}

\subsection{First-stage routing constraints} \label{subsec:2_stage}

The constraints for the first-stage enforce the routing constraints, i.e., the requirements that each POI $i\in T$ should be visited by only one of the UGVs and that each UGV never runs out of fuel as it traverses its route. In the first-stage, the fuel constraint is enforced using the nominal value of fuel consumed by any UGV to traverse any edge $(i,j) \in E$. The first-stage routing constraints are as follows:
\begin{subequations}
	\begin{flalign}
	&x^m(\delta^+(d)) = x^m(\delta^-(d)) \quad \forall d\in D\setminus\{d_0\}, m \in M, \label{f2_1} & \\
	&x^m(\delta^+(d_0)) = 1, \quad \forall m \in M, \label{f2_2} &\\
	& x^m(\delta^-(d_0)) = 1, \quad \forall m \in M, \label{f2_3} &\\
	&x^m(\delta^+(S)) \geqslant 1 \quad \forall S\subset V\setminus\{d_0\} : S\cap T \cap D \neq \emptyset, & \nonumber \\ 
	& \quad \forall m \in M,  \label{f2_4} & \\	
	&\sum_{m \in M} x^m(\delta^+(i)) = 1 \text{ and } \sum_{m \in M} x^m(\delta^-(i)) = 1 \quad \forall  i \in T, \label{f2_5} &\\
	& x^m(\delta^+(i)) = x^m(\delta^-(i))  \quad \forall  i \in T, m \in M, \label{f2_61} &	
\end{flalign}
\label{eq:1stage}
\end{subequations}

\begin{subequations}
	\begin{flalign}	
	& \sum_{j\in V}z_{ij} - \sum_{j\in V}z_{ji}  = \sum_{j\in V}\hat f_{ij}^mx_{ij}^m \quad \forall i \in T, \forall m \in M, \label{f2_6} & \\	
	&z_{di} = \hat f_{di}^mx_{di}^m \quad \forall  i\in V,\, d \in D, m \in M,  \label{f2_7} & 	\\	
	&0 \leqslant z_{ij} \leqslant F_mx_{ij}^m \quad \forall  (i,j) \in E, m \in M,  \label{f2_8} & \\
	&x_{ij}^m \in \{0,1\} \quad \forall  (i,j) \in E, m \in M. \label{f2_9} &
	\end{flalign}
	\label{eq:1stage}
\end{subequations}

Constraint \eqref{f2_1} forces the in-degree and out-degree of each refueling station to be equal. Constraints \eqref{f2_2} and \eqref{f2_3} ensure that all the UGVs leave and return to depot $d_0$. Constraint \eqref{f2_4} ensures that a feasible solution is connected. For each POI $i$, the pair of constraints in \eqref{f2_5} require that some UGV visits the POI $i$. Constraint \eqref{f2_61} forces the in-degree and out-degree of each POI to be equal. Constraint \eqref{f2_6} eliminates subtours of the targets and also defines the flow variables $z_{ij}$ for each edge $(i,j) \in E$ using the nominal fuel consumption values $\hat f_{ij}$. Constraints \eqref{f2_7} -- \eqref{f2_8} together impose the fuel constraints on the routes for all the UGVs. Finally, constraint \eqref{f2_9} imposes binary restrictions on the decision variables $x_{ij}^m$. 

\subsection{Second-stage constraints} \label{f2-ss}
The second-stage model for a fixed $\bm x$ and $\bm f(\omega)$ is given as follows: 
\begin{subequations}
	\begin{flalign}
	&\beta(\bm x,\bm f(\omega)) = \min \sum\limits_{\substack{(i,j) \in E \\ m \in M}}  f_{ij}^m(\omega)\bar{y}_{ij}^m(\omega) \label{ss-f2-1}& \\
	&\text{subject to: } \notag & \\	
	&y^m(\delta^+(d))(\omega) = y^m(\delta^-(d))(\omega) \quad \forall d\in D\setminus\{d_0\}, & \nonumber \\
	& \quad \forall m \in M,  \label{ss-f2-2} & \\		
	&y^m(\delta^+(d_0)) (\omega) = 1, \quad \forall m \in M, \label{ss-f2-3} &\\
	&y^m(\delta^-(d_0)) (\omega) = 1, \quad \forall m \in M, \label{ss-f2-4} &\\
	&y^m(\delta^+(S))(\omega) \geqslant 1 \quad \forall S\subset V\setminus\{d_0\} : S\cap T \cap D \neq \emptyset,  & \nonumber \\
	& \quad \forall m \in M,  \label{ss-f2-5} & \\			
	&\sum_{m \in M} y^m(\delta^+(i)) (\omega) = 1 \text{ and } \sum_{m \in M} y^m(\delta^-(i)) (\omega) = 1,  & \nonumber \\
	& \quad \forall  i \in T,  \label{ss-f2-6} & \\		
	& y^m(\delta^+(i)) = y^m(\delta^-(i))  \quad \forall  i \in T, m \in M, \label{ss-f2-61} &\\		
	& \sum_{j\in V}\bar{z}_{ij}(\omega) - \sum_{j\in V}\bar{z}_{ji}(\omega)  = \sum_{j\in V} f_{ij}^m(\omega)y_{ij}^m (\omega) \quad \forall i \in T, & \nonumber \\
	& \quad \forall m \in M,  \label{ss-f2-7} & 				
	\end{flalign}
\label{eq:1stage1}
\end{subequations}

\begin{subequations}
\begin{flalign}	
	&\bar{z}_{di} = f_{di}^m(\omega) y_{di}^m (\omega) \quad \forall  i\in V,\, d \in D, m \in M,  \label{ss-f2-8} & \\
	&0 \leqslant \bar{z}_{ij}(\omega) \leqslant F_my_{ij}^m(\omega) \quad \forall  (i,j) \in E,  \label{ss-f2-9} & \\
	& \bar{y}_{ij}^m(\omega) \geqslant y_{ij}^m(\omega) - x_{ij}^m \quad \forall (i,j) \in E, m \in M, \label{ss-f2-10} & \\	
	&y_{ij}^m(\omega) \in \{0,1\}, \bar{y}_{ij}^m(\omega) \geqslant 0 \quad \forall  (i,j) \in E \label{ss-f2-11}. &	
	\end{flalign}
	\label{f2-ssc}
\end{subequations}
The decision variables $y_{ij}^m(\omega)$ are similar to the first-stage variable $\bm x$ thus obtaining the paths for each UGV based on the realization $\bm f(\omega)$, and the variables  $\bar{y}_{ij}^m(\omega)$ denote the differences in edges traveled by each UGV `$m$' compared to the first-stage decision $\bm x$. The objective function \eqref{ss-f2-1} minimizes the total swapped and additional travel for the UGVs for a given scenario in comparison to the corresponding first-stage decisions. Constraints \eqref{ss-f2-2}-\eqref{ss-f2-9} are similar to the first-stage constraints. Constraint \eqref{ss-f2-10} estimates the swapped or additional travel for each of the UGVs for the given scenario. Finally, constraint \eqref{ss-f2-11} imposes the restrictions on the second-stage variables.   

\section{Formulation 2} \label{sec:for3} 
\textit{Given a team of UGVs and a subset of it is randomly available for the mission, and a set of POIs sites to visit, find a path for each UGV such that each POI is visited by at most by one UGV, and an objective based on the incentives of POIs visited by the UGVs is maximized}.

\subsection{Objective function} \label{subsec:obj3} The objective function for the two-stage stochastic programming model is the sum of the first-stage profit and the expected second-stage profits. The second-stage profit for a realization $\omega \in \Omega$ of the fuel consumption of the UGV is the change in profit for the realization $\omega$. The reduction in profits is a function of the first-stage routing decision $\bm x$ and the realization $\omega$. Letting the recourse cost be denoted by $\beta(\bm x, \bm z, \bm f(\omega))$, the objective function for the two-stage stochastic optimization problem is given by:

\begin{flalign}
\max \,\, C, \text{ where } C \triangleq \sum\limits_{\substack{(i,j) \in E: j \in T \\ m \in M}} p_{j}^m x_{ij}^m + \mathbb{E}_{\Omega} \left[ \beta(\bm x, \bm z, \bm f) \right] \nonumber \\ =  \sum\limits_{\substack{(i,j) \in E: j \in T \\ m \in M}} p_{j}^m x_{ij}^m + \sum_{\omega \in \Omega} p(\omega) \beta(\bm x, \bm z, \bm f(\omega)). \label{eq:obj}
\end{flalign}

\subsection{First-stage routing constraints} \label{subsec:3_stage}
The constraints for the first-stage enforce the routing constraints, i.e., the requirements that each POI in $T$ can be visited at least once by some UGV and that each UGV never runs out of fuel as it traverses its route. In the first-stage, the fuel constraint is enforced using the nominal value of fuel consumed by any UGV to traverse any edge $(i,j) \in E$. The first-stage routing constraints are as follows:

\begin{subequations}
	\begin{flalign}
	&x^m(\delta^+(d)) = x^m(\delta^-(d)) \quad \forall d\in D\setminus\{d_0\}, m \in M, \label{eq:degree_d_1} & \\
	&x^m(\delta^+(d_0)) = 1, \quad \forall m \in M, \label{eq:degree_d0_1_1} &\\
	& x^m(\delta^-(d_0)) = 1, \quad \forall m \in M, \label{eq:degree_d0_2_1} &\\
	&x^m(S) \leqslant |S|-1 \quad \forall S\subset V\setminus\{d_0\} \nonumber &\\  
	&: S\cap T \cap D \neq \emptyset, \quad \forall m \in M,  \label{eq:sec_1} & \\
	&\sum_{i\in V, m \in M} x_{ij}^m \leqslant 1 \text{ and } \sum_{i\in V, m \in M} x_{ji}^m \leqslant 1 \quad \forall \, j \in T, \label{eq:5} &\\
	& x^m(\delta^+(i)) = x^m(\delta^-(i))  \quad \forall  i \in T, m \in M, \label{eq:51} &\\			
	& \sum_{j\in V}z_{ij}^m - \sum_{j\in V}z_{ji}^m  = \sum_{j\in V}f_{ij}^mx_{ij}^m \quad \forall i \in T, m \in M, \label{eq:7} & \\
	&z_{di}^m = f_{di}^mx_{di}^m \quad \forall \, i\in V, \, d \in D, m \in M, \label{eq:11} & \\	
	&0 \leqslant \sum_{(i,j) \in E} f_{ij}^mx_{ij}^m \leqslant F_m\quad \forall  (i,j) \in E, m \in M, \label{eq:fuel_3_1} & \\
	&x_{ij}^m \in \{0,1\} \, , \,\, z_{ij}^m \leq F_m \quad \forall  (i,j) \in E, m \in M. \label{eq:bin_1}. &
	\end{flalign}
	\label{eq:1stage}
\end{subequations}

Constraint \eqref{eq:degree_d_1} forces the in-degree and out-degree of each refueling station to be equal. Constraints \eqref{eq:degree_d0_1_1} and \eqref{eq:degree_d0_2_1} ensure that all the UGVs leave and return to depot $d_0$, where $m$ is the number of UGVs. Constraint \eqref{eq:sec_1} ensures that a feasible solution is connected. For each target $i$, the pair of constraints in \eqref{eq:5} state that some UGV visits the POI $i$ only once. Constraint \eqref{eq:51} forces the in-degree and out-degree of each POI station to be equal. Constraints \eqref{eq:7}-\eqref{eq:11} eliminates sub-tours of the POIs and also defines the flow variables $z_{ij}^m$ for each edge $(i,j) \in E$ and UV $h$. Constraints \eqref{eq:fuel_3_1} impose the fuel capacity constraints on the routes for all the UGVs. Finally, constraint \eqref{eq:bin_1} imposes restrictions on the decision variables. 

\subsection{Second-stage constraints} \label{subsec:second-stage}
The second-stage model for a fixed $\bm x$, {$\bm z$}, and $\bm f(\omega)$ is given as follows: 
\begin{subequations}
	\begin{flalign}
	&\beta(\bm x, \bm z, \bm f(\omega)) = \max \sum\limits_{\substack{(i,j) \in E \\ m \in M}} -p_{j}^m v_{ij}^m(\omega)  \label{eq:recourse}& \\
	&\text{subject to: } \notag & \\	
	& \sum_{j\in V}f_{ij}^mv_{ij}^m(\omega) = \sum_{j\in V}z_{ij}^m - \sum_{j\in V}z_{ji}^m - \alpha_{m}(\omega)\sum_{j\in V} f_{ij}^m x_{ij}^m   & \nonumber \\
	& \quad \forall  i \in T, m \in M,  \label{eq:14} & \\			
	& v_{ij}^m(\omega) \leq x_{ij}^m \quad \forall (i,j) \in E, m \in M, \label{eq:15a} & \\
	&f_{di}^m v_{di}^m(\omega) = z_{di}^m - \alpha_{m}(\omega).f_{di}^m x_{di}^m \quad \forall \, i\in V, \, d \in D, & \nonumber \\
	& \quad \forall m \in M,  \label{eq:11a} & \\			
	& v_{ij}^m(\omega) \in \{0,1\} \quad \forall \, (i,j) \in E, m \in M.\label{eq:16} & 
	\end{flalign}
	\label{eq:2stage}
\end{subequations}

In the second-stage, $\alpha_{m}(\omega)$ takes a value of $1$ or $0$ denoting the availability of an UGV $m$ for the scenario $\omega$ or not. Variable $v_{ij}^m(\omega)$ maintains the feasibility of the constraints \eqref{eq:14}-\eqref{eq:15a} for the given first-stage values $x$ and $z$. Constraint \eqref{eq:15a} states the dependence of $x_{ij}^m$ and $v_{ij}^m(\omega)$, and finally binary restrictions for $v_{ij}^m(\omega)$ are presented in \eqref{eq:16}. Let the relaxed recourse problem for $\beta(\bm x, \bm z, \bm f(\omega))$ be represented as $\beta_r(\bm x, \bm z, \bm f(\omega))$. In  $\beta_r(\bm x, \bm z, \bm f(\omega))$, the constraints \eqref{eq:16} are replaced by $0 \leq v_{ijh}(\om) \leq 1$.

%

\begin{thm}\label{31}
	The objective values of $\beta(\bm x, \bm z, \bm f(\omega))$ and $\beta_r(\bm x, \bm z, \bm f(\omega))$ are same.
\end{thm}

\section{Algorithm}\label{det-model}

The constraints \eqref{eq:fuel_3_1} are the typical knapsack constraints and the formulation will resemble `orienteering problem'. We will refer the formulation with and without knapsack constraints as TS-OP and TS, respectively. In this section, we present a decomposition algorithm to solve problem TS and TS-OP. The formulations TS and its variants can be provided to any commercial branch-and-cut solvers to obtain an optimal solution. However, observe that the formulations will contain constraint \eqref{eq:sec_1} to ensure any feasible solution is connected. The number of such constraints is exponential and it may not be computationally efficient to enumerate all these constraints and provide them upfront to the solvers. Additionally, stochastic integer programs are  large in scale due to the variables and constraints in the scenarios and they require decomposition algorithms to exploit the special structure of the problem. These challenges and opportunities motivated us to design a decomposition algorithm to solve the instances for TS and its variants. 

\subsection{Decomposition algorithm}\label{det-model}
The decomposition algorithm is a variant of L-shaped algorithm where the deterministic parameters are used to obtain first-stage solutions, and then the second-stage programs are solved based on the obtained first-stage solutions. Then the optimality cuts are generated and added to the first-stage program to approximate the value function of the second-stage cost. The dual information of all the realization of the random data are used to generate the optimality cuts for the first-stage. The use of L-shaped method for TS is possible only due to the theorem \eqref{31}. Otherwise due to the binary restrictions for second-stage variables, the value function will be non-convex and lower semi-continuous in general, and a direct use of L-shaped method is not possible. The first-stage problem is solved as a mixed-integer program with binary restrictions for $x$ variables and by theorem \eqref{31}, the second-stage programs are solved as linear programs. 

\subsubsection{Problem reformulation}\label{decomp-model}	
For the sake of decomposition, the first-stage problem \eqref{eq:obj} -\eqref{eq:bin_1} is reformulated as the following master-problem (MP) and we add an unrestricted variable $\theta$ to the first-stage program. In the formulation TS-MP, let $h(\omega)$ and $\mu(\omega)$ represent the right-hand side and dual values for the second-stage constraints \eqref{eq:14} - \eqref{eq:16} and $v_{ijh}(\om) \in \{0,1\}$ are replaced by $0 \leq v_{ijh}(\om) \leq 1$. Similarly, $T(\omega)$ and $T'(\omega)$ represent the co-efficient matrices for the first-stage variables $x$ and $z$ in the second-stage constraint \eqref{eq:14} - \eqref{eq:16}, respectively. The first-stage master program for the decomposition algorithm TS-MP is given as follows:	

\begin{subequations}
	\begin{flalign}
	& z^ k = \max \sum\limits_{\substack{(i,j) \in E: j \in T \\ m \in M}} p_{j}^m x_{ij}^m + \theta \label{eq:mp} & \\		
	& \text{ Subject to:} \nonumber & \\
	& \eqref{eq:degree_d_1} - \eqref{eq:bin_1}, \label{eq:12} &
	\end{flalign}
\label{eq:1MP}
\end{subequations}

\begin{subequations}
	\begin{flalign}	
	& \sum \limits_{\substack{\omega \in \Omega}} \sum\limits_{\substack{(i,j) \in E, \\ m \in M}} (((\pi_1(\omega)^{t \top}T_1)+(\pi_2(\omega)^{t \top}T_2) \nonumber &\\ 
	& + (\pi_3(\omega)^{t \top}T_3)) x_{ij}^m + ((\pi_1(\omega)^{t \top}S_1)+(\pi_3(\omega)^{t \top}S_3)) z_{ij}^m ) \nonumber &\\ 
	& + \theta \leq \sum \limits_{\substack{\omega \in \Omega}} \pi(\omega)^t h(\omega) , \ t \in 1,...,k \label{eq-master-1a} &\\		
	&x_{ij}^m \in \{0,1\} \, , \,\, z_{ij}^m \leq F_m \quad \forall  (i,j) \in E, m \in M, \,\,\,\theta \in \mathbb{R}. \label{eq:bin_12} &
	\end{flalign}
	\label{eq:1MP}
\end{subequations}

In the master problem \eqref{eq:1MP}, for a scenario $\omega$, $\pi_1(\omega)$, $\pi_2(\omega)$, and $\pi_3(\omega)$ are the dual vectors of the constraints \eqref{eq:14}, \eqref{eq:15a}, and \eqref{eq:11a}, respectively. Similarly, $T_1$, $T_2$, and $T_3$ represent the co-efficient matrices for the variables $x_{ij}^m$ in the constraints \eqref{eq:14}, \eqref{eq:15a}, and \eqref{eq:11a}, respectively. Also, $S_1$ and $S_3$ represent the co-efficient matrices for the variables $z_{ij}^m$ in the constraints \eqref{eq:14}, and \eqref{eq:11a}, respectively. Finally, $\pi(\omega)$ and $h(\omega)$ represent the dual vector and righthand side for the constraints \eqref{eq:14} - \eqref{eq:16}. $\theta$ is an unrestricted decision variable. Constraints \eqref{eq-master-1a} are the \textit{optimality} cuts, which are computed based on the optimal dual solution of all the subproblems given as the second-stage problem $\beta_r(\bm x, \bm z, \bm f(\omega))$. Optimality cuts approximate the value function of the second-stage subproblems $\beta_r(\bm x, \bm z, \bm f(\omega))$. It should be noted that the model TS-MP has relatively complete recourse property, i.e, $\beta_r(\bm x, \bm z, \bm f(\omega)) < \infty$ for any $x_{ij}^m$ and $z_{ij}^m$. 	

We would like to emphasize that the number of optimality cuts generated from second-stage dual values can be a single cut or multi-cut. In single-cut, a cut is generated across all the second-stage problems and in multi-cut, each second-stage program will be approximated by a cut in the first-stage program. In our computational experiments, we adopted a single cut approach as we do not want to stress the first-stage problem as it already consists of binary variables and sub-tour elimination constraints \eqref{eq:sec_1}. Also, the presented algorithm can be extended to instances of TS-OP as the changes occur only in first-stage constraints.

\begin{algorithm}[H]
	\caption{Branch and Cut Algorithm}\label{algo:improvement}
	\begin{algorithmic}[1]
		\doublespacing 
		\State \textbf{Step 0.} Initialize.\\ $n \leftarrow 0$, $lb \leftarrow -\infty$, $\epsilon$ is an user defined parameter and $ub\leftarrow \infty$, and  $x^0, z^0$ are obtained as follows: argmax$ \{ \sum_{(i,j) \in E, m \in M} p_{j}^m x_{ij}^m | x, z \in \eqref{eq:degree_d_1} -\eqref{eq:bin_1}\}.$ \Comment{Initial solution}
		\State \textbf{Step 1.} Solve second-stage programs. \\ Solve $\beta_r(\bm x, \bm z, \bm f(\omega))$ for each $\omega \in \Omega$, and obtain dual solution $\pi(\omega)$ for each second-stage program.  \Comment{Solve sub-problems for given first-stage solution}
		\State \textbf{Step 2.} Optimality cut. 
		Based on the dual values $\pi(\omega)$s from the sub-problems, generate an optimality cut \eqref{eq-master-1a} and add it to the master problem.  \Comment{Sub-problems' objective function are approximated}	
		\State \textbf{Step 3.} Obtain upper bounds.\\
		Solve the master program TS-MP with the new optimality cut  and let the objection function value be $u^n$. Set $ub \leftarrow min\{u^n, ub\}$. \Comment{Sub-problems' objective function are approximated}	
		\State \textbf{Step 4.} Add sub-tour elimination constraints. \\ Check for strongly connected components and if there are any sub-tours, add the constraint \eqref{eq:sec_1}. Go to Step 3. \Comment{Check whether the paths are connected}	
		\State \textbf{Step 5.} Update bounds. \\
		$v^n \leftarrow \{\sum_{(i,j)\in E, m \in M} p_{j}^m x_{ij}^m| x, z \in  \eqref{eq:degree_d_1} -\eqref{eq:bin_1}\} + \mathbb{E}_{\Omega} \left[ \beta(\bm x, \bm z, \bm f) \right]$, and set $lb \leftarrow max\{v^n, lb\}$. If $lb$ is updated, set incumbent solution to $x^{*} \leftarrow x^n,$ and $z^{*}\leftarrow z^n$.
		\Comment{Update the incumbent solution}		
		\State \textbf{Step 6.} Termination. If $ub-lb < \epsilon|ub|$ then stop, $x^{*},$ and $z^{*}$ are the optimal solutions, else set $n \leftarrow n+1$ and return to step 1. \Comment{termination condition}
		\vspace{1ex}
	\end{algorithmic}
\end{algorithm}

\subsection{Sub-tour elimination constraints}\label{det-model}

In our algorithm, we relax the constraints \eqref{eq:sec_1} from the formulation, and whenever the first-stage problem obtains an integer feasible solution to this relaxed problem, we check if any of the constraints \eqref{eq:sec_1} are violated by the integer feasible solution. If so, we add the infeasible constraint to the first-stage problem. This process of adding constraints to the problem sequentially has been observed to be computationally efficient for the TSP, VRP and a huge number of their variants. 

Now, we will detail the algorithm used to find a constraint \eqref{eq:sec_1} that is violated for a given integer feasible solution to the relaxed problem.  A violated constraint \eqref{eq:sec_1} can be described by a subset of vertices $S \subset V\setminus\{d_0\}$ such that $S\cap D \neq \emptyset$ and $x(S) = |S|$ for every $d \in S\cap D$. We find the strongly connected components of $S$. Every strongly connected component that does not contain the depot is a subset $S$ of $V\setminus\{d_0\}$ which violates the constraint \eqref{eq:sec_1}. We add all these infeasible constraints and continue solving the original problem. Many off-the-shelf commercial solvers provide a feature called ``solver callbacks'' to implement such an algorithm into its branch-and-cut framework.


\section{Computational Experiments}\label{det-model}

In this section, we discuss the computational performance of the branch-and-cut algorithm for  formulations presented in the Sec. \ref{sec:formulation}. The mixed-integer linear programs were implemented in Java, using the traditional branch-and-cut framework and the solver callback functionality of CPLEX version 12.6.2. All the simulations were performed on a Dell Precision T5500 workstation (Intel Xeon E5630 processor @2.53 GHz, 12 GB RAM). The computation times reported are expressed in seconds, and we imposed a time limit of 3,600 seconds for each run of the algorithm. The performance of the algorithm was tested with randomly generated test instances. \\

\noindent {\it Instance generation}

The problem instances were randomly generated in a square grid of size [100,100]. The number of refueling stations was set to 4 and the locations of the depot and all the refueling stations were fixed a priori for all the test instances. The number of POIs varies from $10$ to $40$ in steps of five, while their locations were uniformly distributed in the square grid; for each $|T| \in \{10,15,20,25,30,25,40\}$, we generated five random instances. For each of the above generated instances, the number of UGVs in the depot was $3$, and the fuel capacity of the UGVs, $F$, was varied linearly with a parameter $\lambda$. $\lambda$ is defined as the maximum distance between the depot and any POI. The fuel capacity $F$ was assigned a value from the set $\{2.25 \lambda, 2.5 \lambda, 2.75 \lambda, 3\lambda \}$. The travel costs and the fuel consumed to travel between any pair POIs vertices were assumed to be directly proportional to the Euclidean distances between the pair and rounded down to the nearest integer. \\

The utility of the stochastic programming approach can be evaluated by estimating the value of the stochastic solution (VSS) introduced by \cite{birge1982value}. The objective value of the recourse problem (RP) can be stated as in \eqref{eq:obj2}; then we take the expected value of the random variable and solve the \textit{expected value problem} (EV), where $\bm f(\omega)$ in  \eqref{eq:obj2} is replaced by $\bm f(\bar{\omega})$ which represents the mean of the random variable $\bm f(\omega)$. Considering $\bar{x}_{ij}^m$ as the solution for the EV problem, the expected result of using the expected value solutions $(\bar{x})$ is EEV which is given in \eqref{eq:eev}. Then the VSS can be defined as the difference between the objective values of the recourse problem and the EEV, i.e, VSS=EEV-RP.

\begin{flalign}
\min \,\, C, \text{ where } C \triangleq \sum\limits_{\substack{(i,j) \in E \\ m \in M}} \hat f_{ij}^m \bar{x}_{ij}^m + \mathbb{E}_{\Omega} \left[ \beta(\bm \bar{x},\bm f) \right]  \nonumber \\ =  \sum\limits_{\substack{(i,j) \in E \\ m \in M}} \hat f_{ij}^m \bar{x}_{ij}^m + \sum_{\omega \in \Omega} p(\omega) \beta(\bm \bar{x}, \bm f(\omega)). \label{eq:eev}
\end{flalign}

\begin{figure}[!hbtp]
	\begin{center}
		\includegraphics[scale=0.8]{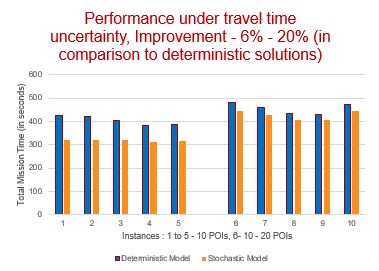}
	\end{center}
	\caption{Performance of VSS while considering uncertainties for travel time among depots and points of interests.}
	\label{Fig2}
\end{figure}

\begin{figure}[!hbtp]
	\begin{center}
		\includegraphics[scale=0.8]{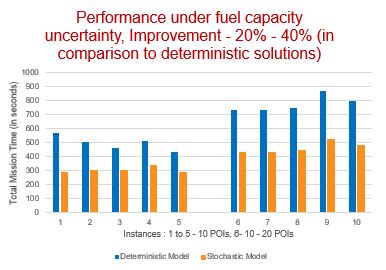}
	\end{center}
	\caption{Performance of VSS while considering uncertainties for fuel capacity.}
	\label{Fig3}
\end{figure}

Figures \ref{Fig2} and \ref{Fig3} represent the results for the formulation presented in section \ref{sec:for2}. In the computational experiments, instances 1 to 5 used 10 POIs and instances 6 to 10 used 20 POIs. Figure \ref{Fig2} shows the use of two-stage model when uncertainties are considered for travel time among depots and points of interests. Gamma distribution is used to characterize the uncertainty in travel time and continuous beta distribution for the uncertainty in fuel consumption. The results from stochastic model are compared with deterministic model. Overall, under travel time uncertainty, the average improvement indicated by VSS is between 6\% and 20\%. 

Figure \ref{Fig3} demonstrates the use of two-stage model when uncertainties are considered for fuel capacity. Gamma distribution is used to characterize uncertainty in travel time and continuous beta distribution for uncertainty in fuel consumption. The results from stochastic model are compared with deterministic model. Overall, under fuel capacity uncertainty, the average improvement indicated by VSS is between 20\% and 40\%.

\begin{figure}[!hbtp]
	\begin{center}
		\includegraphics[scale=0.5]{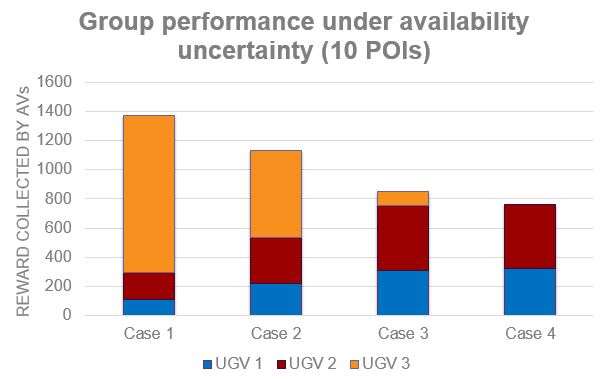}
	\end{center}
	\caption{ Uncertainty in availability of UGVs - 10 POIs}
	\label{Fig4}
\end{figure}

\begin{figure}[!hbtp]
	\begin{center}
		\includegraphics[scale=0.5]{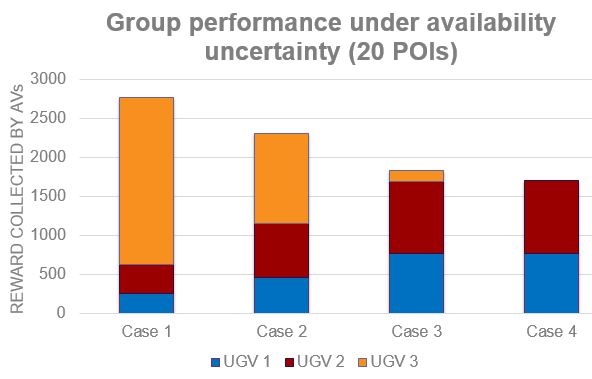}
	\end{center}
	\caption{ Uncertainty in availability of UGVs - 20 POIs.}
	\label{Fig5}
\end{figure}

Figures \ref{Fig4} and \ref{Fig5} represent the results for the formulation presented in the section \ref{sec:for3}. In the computational experiments,  figures \ref{Fig4} and \ref{Fig5} used 10 and 20 POIs, respectively. Each POI has a reward and is visited by at most one UGV, and the objective is to maximize the total reward collected by the UGVs. UGV 3’s availability with probabilities are 1, 0.75, 0.25, and 0 in Cases 1, 2, 3, and 4, respectively. As represented in the figures, the contribution of UGV 3 monotonically decreases depending upon the  probability of availability. This type of marginal decrease in the contribution of  UGV 3 is not possible with a deterministic model since the availability has to be considered as binary (0 or 1). Hence, the deterministic model can only handle extreme cases, and a stochastic model can include marginal increase or decrease of an asset's utilization based on the probability distribution to denote its availability.

\section{Conclusion}

Path planning problem for manned and unmanned UGVs is an important area of research for efficient use of UGVs. This paper presents two different stochastic programming models to address uncertainties; travel time and availability of UGVs. To overcome the computational complexity, a decomposition algorithm is presented. Computational experiments are performed to demonstrate the usefulness of stochastic models over their deterministic versions. A potential future study is to evaluate the robustness of the solutions under cost minimization and profit maximization, and choose an appropriate objective based on the uncertainties in the environment. 

\section{Acknowledgement}
The authors wish to acknowledge the technical and financial support of the Automotive Research Center (ARC) in accordance with Cooperative Agreement W56HZV-19-2-0001 U.S. Army CCDC Ground Vehicle Systems Center (GVSC) Warren, MI. 

\section{Legal Statement}
DISTRIBUTION  A.  Approved  for  public  release:  distribution  unlimited.

\section{References}

\bibliography{Referencias}


\end{document}